\documentclass[11pt]{amsart}
\usepackage{amsmath,amssymb,mathrsfs}
\newtheorem{theorem}{Theorem}
\newtheorem{proposition}[theorem]{Proposition}
\newtheorem{corollary}[theorem]{Corollary}
\newtheorem{lemma}[theorem]{Lemma}
\newtheorem*{minoo}{Min-Oo's Conjecture}

\begin{document}

\title[Deformations that increase scalar curvature]{Deformations of the hemisphere that increase scalar curvature}
\author{Simon Brendle, Fernando C. Marques, and Andre Neves}
\address{Department of Mathematics \\ Stanford University \\ 450 Serra Mall, Bldg 380 \\ Stanford, CA 94305} 
\address{Instituto de Matem\'atica Pura e Aplicada (IMPA) \\ Estrada Dona Castorina 110 \\ 22460-320 Rio de Janeiro \\ Brazil}
\address{Imperial College \\ Huxley Building \\ 180 Queen's Gate \\ London SW7 2RH \\ United Kingdom}
\thanks{The first author was supported in part by the National Science Foundation under grant DMS-0905628. The second author was supported by CNPq-Brazil, FAPERJ, and the Stanford Department of Mathematics.}
\begin{abstract}
Consider a compact Riemannian manifold $M$ of dimension $n$ whose boundary $\partial M$ is totally geodesic and is isometric to the standard sphere $S^{n-1}$. A natural conjecture of Min-Oo asserts that if the scalar curvature of $M$ is at least $n(n-1)$, then $M$ is isometric to the hemisphere $S_+^n$ equipped with its standard metric. This conjecture is inspired by the positive mass theorem in general relativity, and has been verified in many special cases.

In this paper, we construct counterexamples to Min-Oo's Conjecture in dimension $n \geq 3$.
\end{abstract}
\maketitle

\section{Introduction}

One of the major results in differential geometry is the positive mass theorem, which asserts that any asymptotically flat manifold $M$ of dimension $n \leq 7$ with nonnegative scalar curvature has nonnegative ADM mass. Furthermore, the ADM mass is strictly positive unless $M$ is isometric to the Euclidean space $\mathbb{R}^n$. This theorem was proved in 1979 by Schoen and Yau \cite{Schoen-Yau1} using minimal surface techniques. Witten \cite{Witten} subsequently gave an alternative proof of the positive mass theorem based on spinors and the Dirac operator (see also \cite{Bartnik1}, \cite{Parker-Taubes}). Witten's argument works for any spin manifold $M$, without any restriction on the dimension. Similar techniques can be used to show that the torus $T^n$ does not admit a metric of positive scalar curvature (see \cite{Gromov-Lawson1}, \cite{Gromov-Lawson2}, \cite{Schoen-Yau2}, \cite{Schoen-Yau3}).

It was observed by Miao \cite{Miao} that the positive mass theorem implies the following rigidity result for metrics on the unit ball:

\begin{theorem}
\label{rigidity.for.ball}
Suppose that $g$ is a smooth metric on the unit ball $B^n \subset \mathbb{R}^n$ with the following properties: 
\begin{itemize}
\item The scalar curvature of $g$ is nonnegative.
\item The induced metric on the boundary $\partial B^n$ agrees with the standard metric on $\partial B^n$.
\item The mean curvature of $\partial B^n$ with respect to $g$ is at least $n-1$.
\end{itemize}
Then $g$ is isometric to the standard metric on $B^n$.
\end{theorem}

Theorem \ref{rigidity.for.ball} was generalized by Shi and Tam \cite{Shi-Tam}. The following result is a subcase of Shi and Tam's theorem (cf. \cite{Shi-Tam}, Theorem 4.1):

\begin{theorem}[Y.~Shi, L.F.~Tam]
\label{shitam}
Let $\Omega$ be a strictly convex domain in $\mathbb{R}^n$ with smooth boundary. Moreover, suppose that $g$ is a Riemannian metric on $\Omega$ with the following properties: 
\begin{itemize}
\item The scalar curvature of $g$ is nonnegative.
\item The induced metric on the boundary $\partial \Omega$ agrees with the restriction of the Euclidean metric to $\partial \Omega$.
\item The mean curvature of $\partial \Omega$ with respect to $g$ is positive.
\end{itemize}
Then 
\[\int_{\partial \Omega} (H_0 - H_g) \, d\sigma_g \geq 0,\] 
where $H_g$ denotes the mean curvature of $\partial \Omega$ with respect to $g$ and $H_0$ denotes the mean curvature of $\partial \Omega$ with respect to the Euclidean metric. Finally, if equality holds, then $g$ is isometric to the Euclidean metric.
\end{theorem}

Similar rigidity results are known for asymptotically hyperbolic manifolds with scalar curvature at least $-n(n-1)$. The first result in this direction was obtained by Min-Oo \cite{Min-Oo1} in 1989. This result was subsequently extended by Andersson and Dahl \cite{Andersson-Dahl}. There also is an analogue of the positive mass theorem for asymptotically hyperbolic manifolds, due to Chru\'sciel and Herzlich \cite{Chrusciel-Herzlich}, Chru\'sciel and Nagy \cite{Chrusciel-Nagy}, and Wang \cite{Wang}. Finally, Boualem and Herzlich have proved a scalar curvature rigidity theorem for K\"ahler manifolds that are asymptotic to complex hyperbolic space (cf. \cite{Boualem-Herzlich}, \cite{Herzlich}). 

Motivated by the positive mass theorem and its analogue in the asymptotically hyperbolic setting (cf. \cite{Min-Oo1}), Min-Oo proposed the following conjecture (cf. \cite{Min-Oo2}, Theorem 4):

\begin{minoo}
Suppose that $g$ is a smooth metric on the hemisphere $S_+^n$ with the following properties: 
\begin{itemize}
\item The scalar curvature of $g$ is at least $n(n-1)$.
\item The induced metric on the boundary $\partial S_+^n$ agrees with the standard metric on $\partial S_+^n$.
\item The boundary $\partial S_+^n$ is totally geodesic with respect to $g$.
\end{itemize}
Then $g$ is isometric to the standard metric on $S_+^n$.
\end{minoo}

Min-Oo's Conjecture is very natural given the analogy with the positive mass theorem, and was widely expected to be true; see e.g. \cite{Gromov}, p.~47, or \cite{Hang-Wang2}, p.~629. Various attempts have been made to prove Min-Oo's conjecture (using both spinor and minimal surface techniques), and many partial results have been obtained. It follows from a theorem of Toponogov \cite{Toponogov} that Min-Oo's Conjecture holds in dimension $2$ (see also \cite{Hang-Wang2}, Theorem 4). For $n \geq 3$, Hang and Wang \cite{Hang-Wang1} proved that Min-Oo's Conjecture holds for any metric $g$ which is conformally equivalent to the standard metric on $S_+^n$. They also noted that the rigidity statement fails if the hemisphere is replaced by a geodesic ball in $S^n$ of radius strictly greater than $\frac{\pi}{2}$. Huang and Wu \cite{Huang-Wu} showed that Min-Oo's Conjecture holds for a class of hypersurfaces in $\mathbb{R}^{n+1}$ which includes graphs. Eichmair \cite{Eichmair} confirmed the conjecture for $n=3$, assuming that the boundary satisfies an isoperimetric condition. The argument in \cite{Eichmair} uses isoperimetric surfaces; this technique was originally developed by Bray \cite{Bray1} to prove a volume comparison theorem involving scalar curvature (see also \cite{Brendle-survey}). Llarull \cite{Llarull} has established an interesting rigidity result for metrics on $S^n$ with scalar curvature $R_g \geq n(n-1)$ (see also \cite{Gromov}). Llarull's theorem was recently generalized by Listing \cite{Listing}. Finally, the first two authors recently obtained scalar curvature rigidity results for certain geodesic balls in $S^n$; see \cite{Brendle-Marques2} for details.

In 2009, Hang and Wang \cite{Hang-Wang2} proved the following beautiful rigidity theorem:

\begin{theorem}[F.~Hang, X.~Wang]
\label{hangwang.thm}
Suppose that $g$ is a smooth metric on the hemisphere $S_+^n \subset S^n$ with the following properties: 
\begin{itemize}
\item The Ricci curvature of $g$ is bounded from below by $\text{\rm Ric}_g \geq (n-1) \, g$.
\item The induced metric on the boundary $\partial S_+^n$ agrees with the standard metric on $\partial S_+^n$.
\item The second fundamental form of the boundary $\partial S_+^n$ with respect to $g$ is nonnegative.
\end{itemize}
Then $g$ is isometric to the standard metric on $S_+^n$.
\end{theorem}

The proof of Theorem \ref{hangwang.thm} relies on a very interesting application of Reilly's formula (see e.g. \cite{Li}, Section 8). 

Theorem \ref{hangwang.thm} is similar to Min-Oo's Conjecture, except that the lower bound for the scalar curvature is replaced by a lower bound for the Ricci tensor.

In this paper, we construct counterexamples to Min-Oo's Conjecture for each $n \geq 3$. To that end, we proceed in two steps. In a first step, we show that the standard metric on the hemisphere $S_+^n$ can be perturbed so that the scalar curvature increases and the mean curvature of the boundary becomes positive:

\begin{theorem}
\label{thm.A}
Given any integer $n \geq 3$, there exists a smooth metric $g$ on the hemisphere $S_+^n$ with the following properties:
\begin{itemize}
\item The scalar curvature of $g$ is strictly greater than $n(n-1)$.
\item At each point on $\partial S_+^n$, we have $g - \overline{g} = 0$, where $\overline{g}$ denotes the standard metric on $S_+^n$.
\item The mean curvature of $\partial S_+^n$ with respect to $g$ is strictly positive (i.e. the mean curvature vector is inward-pointing).
\end{itemize}
\end{theorem}

The proof of Theorem \ref{thm.A} relies on a perturbation analysis, which is reminiscent of the construction of counterexamples to Schoen's Compactness Conjecture for the Yamabe problem (cf. \cite{Brendle}, \cite{Brendle-Marques1}). 

Lohkamp \cite{Lohkamp2} proved that for any Riemannian manifold there exist local deformations of the metric which decrease scalar curvature (see also \cite{Lohkamp1}). By contrast, it is not always possible to increase scalar curvature. In our case, the main difficulty is that the hemisphere $(S_+^n,\overline{g})$ is static (cf. \cite{Bartnik2}, \cite{Corvino}). As a result of that, the linearization of the scalar curvature fails to be surjective. In particular, Corvino's theorem concerning local deformations of the scalar curvature (see \cite{Corvino}, Theorem 1) does not apply in this situation.

In a second step, we perform another perturbation to make the boundary totally geodesic. Since the examples constructed in Theorem \ref{thm.A} have positive mean curvature, this deformation can be done in such a way that the scalar curvature remains greater than $n(n-1)$. More precisely, we prove the following general result:

\begin{theorem}
\label{thm.B}
Let $M$ be a compact manifold of dimension $n$ with boundary $\partial M$, and let $g$ and $\tilde{g}$ be two smooth Riemannian metrics on $M$ such that $g - \tilde{g} = 0$ at each point on $\partial M$. Moreover, we assume that $H_g - H_{\tilde{g}} > 0$ at each point on $\partial M$. Given any real number $\varepsilon > 0$ and any neighborhood $U$ of $\partial M$, there exists a smooth metric $\hat{g}$ on $M$ with the following properties: 
\begin{itemize}
\item We have the pointwise inequality $R_{\hat{g}}(x) \geq \min \{R_g(x),R_{\tilde{g}}(x)\} - \varepsilon$ at each point $x \in M$.
\item $\hat{g}$ agrees with $g$ outside $U$.
\item $\hat{g}$ agrees with $\tilde{g}$ in a neighborhood of $\partial M$.
\end{itemize}
\end{theorem}

The proof of Theorem \ref{thm.B} involves a purely local construction based on cut-off functions. We expect that Theorem \ref{thm.B} will be useful in other settings.

We note that Bray \cite{Bray2} and Miao \cite{Miao} have used different methods to smooth out a Riemannian metric admitting corners along a hypersurface. However, the construction in \cite{Miao} does not maintain positive scalar curvature; indeed, the scalar curvature of the mollified metric may be negative in a small region (see \cite{Miao}, Proposition 3.1, for details).

Combining Theorem \ref{thm.A} with Theorem \ref{thm.B}, we obtain counterexamples to Min-Oo's Conjecture in dimension $n \geq 3$:

\begin{corollary}
\label{counterexample}
Given any integer $n \geq 3$, there exists a smooth metric $\hat{g}$ on the hemisphere $S_+^n$ with the following properties:
\begin{itemize}
\item The scalar curvature of $\hat{g}$ is strictly greater than $n(n-1)$.
\item At each point on $\partial S_+^n$, we have $\hat{g} - \overline{g} = 0$, where $\overline{g}$ denotes the standard metric on $S_+^n$.
\item The boundary $\partial S_+^n$ is totally geodesic with respect to $\hat{g}$.
\end{itemize}
Moreover, the metric $\hat{g}$ can be chosen to be rotationally symmetric in a neighborhood of $\partial S_+^n$.
\end{corollary}

Let us briefly describe how Corollary \ref{counterexample} follows from Theorem \ref{thm.A} and Theorem \ref{thm.B}. Let $g$ be the metric constructed in Theorem \ref{thm.A}. It is not difficult to construct a rotationally symmetric metric $\tilde{g}$ on $S_+^n$ with the following properties:
\begin{itemize}
\item The scalar curvature of $\tilde{g}$ is strictly greater than $n(n-1)$ in a neighborhood of $\partial S_+^n$.
\item We have $\tilde{g} - \overline{g} = 0$ at each point on $\partial S_+^n$.
\item The boundary $\partial S_+^n$ is totally geodesic with respect to $\tilde{g}$.
\end{itemize}
Applying Theorem \ref{thm.B}, we obtain a metric $\hat{g}$ with the required properties.

Using a slightly different choice of $\tilde{g}$, we can arrange for the boundary $\partial S_+^n$ to be strictly convex (instead of totally geodesic).

We next construct metrics on the hemisphere $S_+^n$ which have scalar curvature at least $n(n-1)$ and agree with the standard metric in a neighborhood of the equator.

\begin{theorem}
\label{thm.C}
Given any integer $n \geq 3$, there exists a smooth metric $\hat{g}$ on the hemisphere $S_+^n$ with the following properties:
\begin{itemize}
\item The scalar curvature of $\hat{g}$ is at least $n(n-1)$ at each point on $S_+^n$.
\item The scalar curvature of $\hat{g}$ is strictly greater than $n(n-1)$ at some point on $S_+^n$.
\item The metric $\hat{g}$ agrees with the standard metric $\overline{g}$ in a neighborhood of $\partial S_+^n$.
\end{itemize}
\end{theorem}

These examples show that there is no analogue of the positive mass theorem in the spherical setting.

Using Theorem \ref{thm.C} and a doubling argument, we obtain the following result:

\begin{corollary}
\label{projective.space}
Given any integer $n \geq 3$, there exists a smooth metric $g$ on the real projective space $\mathbb{RP}^n$ with the following properties:
\begin{itemize}
\item The scalar curvature of $g$ is at least $n(n-1)$ at each point on $\mathbb{RP}^n$.
\item The scalar curvature of $g$ is strictly greater than $n(n-1)$ at some point on $\mathbb{RP}^n$.
\item The metric $g$ agrees with the standard metric in a neighborhood of the equator in $\mathbb{RP}^n$.
\end{itemize}
\end{corollary}

Corollary \ref{projective.space} is of interest in light of recent work of Bray, Brendle, Eichmair, and Neves \cite{Bray-Brendle-Eichmair-Neves} (see also \cite{Bray-Brendle-Neves}). To describe this result, consider a Riemannian metric on $\mathbb{RP}^3$ with scalar curvature at least $6$. The main result of \cite{Bray-Brendle-Eichmair-Neves} asserts that an area-minimizing surface homeomorphic to $\mathbb{RP}^2$ has area at most $2\pi$. Furthermore, equality holds if and only if the metric on $\mathbb{RP}^3$ has constant sectional curvature $1$. Corollary \ref{projective.space} shows that the area-minimizing condition in \cite{Bray-Brendle-Eichmair-Neves} cannot be replaced by stability.

\section{Proof of Theorem \ref{thm.A}}

Let $S^n$ denote the unit sphere in $\mathbb{R}^{n+1}$, and let $\overline{g}$ be the standard metric on $S^n$. Moreover, let $f: S^n \to \mathbb{R}$ denote the restriction of the coordinate function $x_{n+1}$ to $S^n$. For abbreviation, we denote by $S_+^n = \{f \geq 0\}$ the upper hemisphere, and by $\Sigma = \{f = 0\}$ the equator in $S^n$. 

For any Riemannian metric $g$ on $S_+^n$, we denote by $R_g$ the scalar curvature of $g$. Moreover, we denote by $H_g$ the mean curvature of $\Sigma$ with respect to $g$. In other words, the mean curvature vector of $\Sigma$ is given by $-H_g \, \nu_g$, where $\nu_g$ denotes the outward-pointing normal vector to $\Sigma$.

We next define a functional $\mathscr{F}$ on the space of Riemannian metrics by 
\[\mathscr{F}(g) = \int_{S_+^n} R_g \, f \, d\text{\rm vol}_{\overline{g}} + 2 \, \text{\rm area}(\Sigma,g).\] 
We note that similar ideas were used in the work of Fischer and Marsden \cite{Fischer-Marsden}, where the case of manifolds without boundary is studied.

\begin{proposition}
\label{first.variation}
The first variation of the functional $\mathscr{F}$ at $\overline{g}$ vanishes. In other words, if $g(t)$ is a smooth one-parameter family of Riemannian metrics on $S_+^n$ with $g(0) = \overline{g}$, then $\frac{d}{dt} \mathscr{F}(g(t)) \big |_{t=0} = 0$.
\end{proposition}

\textbf{Proof.} 
Let $h = \frac{\partial}{\partial t} g(t) \big |_{t=0}$. Using Theorem 1.174 in \cite{Besse}, we obtain  
\[\frac{\partial}{\partial t} R_{g(t)} \Big |_{t=0} = \sum_{i,j=1}^n (\overline{D}_{e_i,e_j}^2 h)(e_i,e_j) - \Delta_{\overline{g}} (\text{\rm tr}_{\overline{g}}(h)) - (n-1) \, \text{\rm tr}_{\overline{g}}(h).\] 
Here, $\overline{D}$ denotes the Levi-Civita connection with respect to the metric $\overline{g}$, and $\{e_1,\hdots,e_n\}$ is a local orthonormal frame with respect to $\overline{g}$. This implies 
\begin{align*} 
&\frac{d}{dt} \bigg ( \int_{S_+^n} R_{g(t)} \, f \, d\text{\rm vol}_{\overline{g}} \bigg ) \bigg |_{t=0} \\ 
&= \int_{S_+^n} \sum_{i,j=1}^n (\overline{D}_{e_i,e_j}^2 h)(e_i,e_j) \, f \, d\text{\rm vol}_{\overline{g}} \\ &- \int_{S_+^n} \Delta_{\overline{g}}(\text{\rm tr}_{\overline{g}}(h)) \, f \, d\text{\rm vol}_{\overline{g}} - (n-1) \int_{S_+^n} \text{\rm tr}_{\overline{g}}(h) \, f \, d\text{\rm vol}_{\overline{g}} \\ 
&= \int_{S_+^n} \langle h,\overline{D}^2 f \rangle \, d\text{\rm vol}_{\overline{g}} - \int_{S_+^n} \text{\rm tr}_{\overline{g}}(h) \, \Delta_{\overline{g}} f \, d\text{\rm vol}_{\overline{g}} - (n-1) \int_{S_+^n} \text{\rm tr}_{\overline{g}}(h) \, f \, d\text{\rm vol}_{\overline{g}} \\ 
&- \int_\Sigma h(\nabla f,\nu) \, d\sigma_{\overline{g}} + \int_\Sigma \text{\rm tr}_{\overline{g}}(h) \, \langle \nabla f,\nu \rangle \, d\sigma_{\overline{g}}, 
\end{align*}
where $\nu$ denotes the outward-pointing unit normal vector with respect to $\overline{g}$. Clearly, $\nu = -\nabla f$. Using the identity $\overline{D}^2 f = -f \, \overline{g}$, we obtain 
\[\overline{D}^2 f - (\Delta_{\overline{g}} f) \, \overline{g} - (n-1) \, f \, \overline{g} = 0,\] 
hence 
\[\langle h,\overline{D}^2 f \rangle - \text{\rm tr}_{\overline{g}}(h) \, \Delta_{\overline{g}} f - (n-1) \, \text{\rm tr}_{\overline{g}}(h) \, f = 0.\] 
Putting these facts together, we conclude that 
\[\frac{d}{dt} \bigg ( \int_{S_+^n} R_{g(t)} \, f \, d\text{\rm vol}_{\overline{g}} \bigg ) \bigg |_{t=0} = -\int_\Sigma (\text{\rm tr}_{\overline{g}}(h) - h(\nu,\nu)) \, d\sigma_{\overline{g}}.\] 
On the other hand, 
\[\frac{d}{dt} \text{\rm area}(\Sigma,g(t)) \Big |_{t=0} = \frac{1}{2} \int_\Sigma \text{\rm tr}_{\overline{g}}(h|_\Sigma) \, d\sigma_{\overline{g}} = \frac{1}{2} \int_\Sigma (\text{\rm tr}_{\overline{g}}(h) - h(\nu,\nu)) \, d\sigma_{\overline{g}}.\] 
Putting these facts together, we obtain 
\[\frac{d}{dt} \bigg ( \int_{S_+^n} R_{g(t)} \, f \, d\text{\rm vol}_{\overline{g}} + 2 \, \text{\rm area}(\Sigma,g(t)) \bigg ) \bigg |_{t=0} = 0,\] 
as claimed. \\

The following proposition is one of the key geometric ingredients in the argument. It implies that, for $n \geq 3$, there exist deformations of the equator in $S^n$ which increase area and have positive mean curvature. This is the only point in our construction where the condition $n \geq 3$ is used.

\begin{proposition}
\label{eta}
Assume that $n \geq 3$. Then there exists a function $\eta: \Sigma \to \mathbb{R}$ such that 
\[\Delta_\Sigma \eta + (n-1) \eta < 0\] 
and 
\[\int_\Sigma (|\nabla_\Sigma \eta|^2 - (n-1) \eta^2) \, d\sigma_{\overline{g}} > 0.\] 
\end{proposition}

\textbf{Proof.} 
We define a function $\psi: \Sigma \to \mathbb{R}$ by 
\[\psi = -1 + \frac{n-1}{2} \, x_n^2 + \frac{(n-1)(n+1)}{24} \, x_n^4 + \frac{(n-1)(n+1)(n+3)}{240} \, x_n^6.\] 
Using the identities 
\[\Delta_\Sigma x_n + (n-1) x_n = 0\] 
and 
\[x_n^2 + |\nabla_\Sigma x_n|^2 = 1,\] 
we obtain 
\[\Delta_\Sigma \psi + (n-1) \psi = -\frac{(n-1)(n+1)(n+3)(n+5)}{48} \, x_n^6 \leq 0.\] 
We next show that 
\[\int_\Sigma (|\nabla_\Sigma \psi|^2 - (n-1) \psi^2) \, d\sigma_{\overline{g}} > 0.\] 
To that end, we use the following recursive relation: 
\begin{align*} 
\int_\Sigma x_n^\alpha \, d\sigma_{\overline{g}} 
&= \int_\Sigma x_n^\alpha \, (x_n^2 + |\nabla_\Sigma x_n|^2) \, d\sigma_{\overline{g}} \\ 
&= \int_\Sigma x_n^{\alpha+2} \, d\sigma_{\overline{g}} + \frac{1}{\alpha+1} \int_\Sigma \langle \nabla_\Sigma (x_n^{\alpha+1}),\nabla_\Sigma x_n \rangle \, d\sigma_{\overline{g}} \\ 
&= \int_\Sigma x_n^{\alpha+2} \, d\sigma_{\overline{g}} - \frac{1}{\alpha+1} \int_\Sigma x_n^{\alpha+1} \, \Delta_\Sigma x_n \, d\sigma_{\overline{g}} \\ 
&= \frac{n+\alpha}{\alpha+1} \int_\Sigma x_n^{\alpha+2} \, d\sigma_{\overline{g}}. 
\end{align*} 
This implies 
\begin{align*} 
&\int_\Sigma x_n^6 \, d\sigma_{\overline{g}} = \frac{n+6}{7} \int_\Sigma x_n^8 \, d\sigma_{\overline{g}} \\ 
&\int_\Sigma x_n^{10} \, d\sigma_{\overline{g}} = \frac{9}{n+8} \int_\Sigma x_n^8 \, d\sigma_{\overline{g}} \\ 
&\int_\Sigma x_n^{12} \, d\sigma_{\overline{g}} = \frac{99}{(n+8)(n+10)} \int_\Sigma x_n^8 \, d\sigma_{\overline{g}}. 
\end{align*} From this, we deduce that 
\begin{align*} 
\int_\Sigma \psi \, x_n^6 \, d\sigma_{\overline{g}} 
&= \bigg ( -\frac{n+6}{7} + \frac{n-1}{2} + \frac{3(n-1)(n+1)}{8(n+8)} \\ 
&\hspace{10mm} + \frac{33(n-1)(n+1)(n+3)}{80(n+8)(n+10)} \bigg ) \, \int_\Sigma x_n^8 \, d\sigma_{\overline{g}} > 0 
\end{align*}
for $n \geq 3$. Thus, we conclude that 
\begin{align*} 
&\int_\Sigma (|\nabla_\Sigma \psi|^2 - (n-1) \psi^2) \, d\sigma_{\overline{g}} \\ 
&= -\int_\Sigma \psi \, (\Delta_\Sigma \psi + (n-1)\psi) \, d\sigma_{\overline{g}} \\ 
&= \frac{(n-1)(n+1)(n+3)(n+5)}{48} \int_\Sigma \psi \, x_n^6 \, d\sigma_{\overline{g}} > 0 
\end{align*} 
for $n \geq 3$. Hence, if $c > 0$ is sufficiently small, then the function $\eta = \psi - c$ has the required properties. This completes the proof. \\

In the remainder of this section, we will always assume that $n \geq 3$. Let $\eta: \Sigma \to \mathbb{R}$ be the function constructed in Proposition \ref{eta}. Moreover, let $X$ be a smooth vector field on $S^n$ such that 
\[X = \eta \, \nu\] 
and 
\[\overline{D}_\nu X = -\nabla_\Sigma \eta\] 
at each point on $\Sigma$. The vector field $X$ generates a one-parameter group of diffeomorphisms, which we denote by $\varphi_t: S^n \to S^n$. For each $t$, we define two Riemannian metrics $g_0(t)$ and $g_1(t)$ by 
\[g_0(t) = \overline{g} + t \, \mathscr{L}_X \overline{g}\] 
and 
\[g_1(t) = \varphi_t^*(\overline{g}).\] 
It follows from our choice of $X$ that $\mathscr{L}_X \overline{g} = 0$ at each point on $\Sigma$. This implies $g_0(t) - \overline{g} = 0$ at each point on $\Sigma$. Moreover, we have $R_{g_1(t)} = n(n-1)$ for all $t$.

\begin{proposition}
\label{key}
We have 
\[\frac{\partial}{\partial t} R_{g_0(t)} \Big |_{t=0} = 0.\] 
Moreover, the function 
\[Q = \frac{\partial^2}{\partial t^2} R_{g_0(t)} \Big |_{t=0}\] 
satisfies $\int_{S_+^n} Q \, f \, d\text{\rm vol}_{\overline{g}} > 0$. 
\end{proposition}

\textbf{Proof.} 
Clearly, 
\[g_0(0) = g_1(0) = \overline{g}\] 
and 
\[\frac{\partial}{\partial t} g_0(t) \Big |_{t=0} = \frac{\partial}{\partial t} g_1(t) \Big |_{t=0} = \mathscr{L}_X \overline{g}.\] 
This implies 
\[\frac{\partial}{\partial t} R_{g_0(t)} \Big |_{t=0} = \frac{\partial}{\partial t} R_{g_1(t)} \Big |_{t=0} = 0.\] 
This proves the first statement.

We now describe the proof of the second statement. By Proposition \ref{first.variation}, the first variation of the functional $\mathscr{F}$ at $\overline{g}$ vanishes. This implies 
\[\frac{d^2}{dt^2} \mathscr{F}(g_0(t)) \Big |_{t=0} = \frac{d^2}{dt^2} \mathscr{F}(g_1(t)) \Big |_{t=0}.\] 
Since $g_0(t)$ agrees with the standard metric $\overline{g}$ at each point on $\Sigma$, we obtain 
\[\mathscr{F}(g_0(t)) = \int_{S_+^n} R_{g_0(t)} \, f \, d\text{\rm vol}_{\overline{g}} + 2 \, \text{\rm area}(\Sigma,\overline{g}),\] 
hence 
\[\frac{d^2}{dt^2} \mathscr{F}(g_0(t)) \Big |_{t=0} = \int_{S_+^n} Q \, f \, d\text{\rm vol}_{\overline{g}}.\] 
On the other hand, the identity $R_{g_1(t)} = n(n-1)$ implies 
\[\mathscr{F}(g_1(t)) = \int_{S_+^n} n(n-1) \, f \, d\text{\rm vol}_{\overline{g}} + 2 \, \text{\rm area}(\varphi_t(\Sigma),\overline{g}).\] 
Using the standard formula for the second variation of area (see e.g. \cite{Li}, Section 1), we obtain 
\begin{align*} 
\frac{d^2}{dt^2} \mathscr{F}(g_1(t)) \Big |_{t=0} 
&= 2 \, \frac{d^2}{dt^2} \text{\rm area}(\varphi_t(\Sigma),\overline{g}) \Big |_{t=0} \\ 
&= 2 \int_\Sigma (|\nabla_\Sigma \eta|^2 - (n-1) \eta^2) \, d\sigma_{\overline{g}} > 0. 
\end{align*} Putting these facts together, the assertion follows. \\

For abbreviation, we define 
\[\mu = \frac{\int_{S_+^n} Q \, f \, d\text{\rm vol}_{\overline{g}}}{\int_{S_+^n} f \, d\text{\rm vol}_{\overline{g}}}.\] 
It follows from Proposition \ref{key} that $\mu$ is positive. Moreover, we have 
\[\int_{S_+^n} (Q - \mu) \, f \, d\text{\rm vol}_{\overline{g}} = 0\] 
by definition of $\mu$.

\begin{proposition} 
There exists a smooth function $u: S_+^n \to \mathbb{R}$ such that 
\[\Delta_{\overline{g}} u + nu = Q - \mu\] 
and $u|_\Sigma = 0$.
\end{proposition}

\textbf{Proof.}
Let $L$ denote the operator $\Delta_{\overline{g}} + n$ on the hemisphere $S_+^n$ with Dirichlet boundary condition. It is easy to see that $L$ is symmetric and $f$ lies in the nullspace of $L$. Since $f$ is a positive function, we conclude that the first eigenvalue of $L$ is equal to $0$. Therefore, the nullspace of $L$ is one-dimensional, and is spanned by the function $f$. From this, the assertion follows. \\

We now define 
\[g(t) = \overline{g} + t \, \mathscr{L}_X \overline{g} + \frac{1}{2(n-1)} \, t^2 \, u \, \overline{g}.\] 
It is straightforward to verify that $g(t) - \overline{g} = 0$ at each point on $\Sigma$.

The following result implies that the scalar curvature of $g(t)$ is greater than $n(n-1)$ if $t > 0$ is sufficiently small:

\begin{proposition}
We have $\frac{\partial}{\partial t} R_{g(t)} \big |_{t=0} = 0$ and $\frac{\partial^2}{\partial t^2} R_{g(t)} \big |_{t=0} = \mu > 0$ at each point on $S_+^n$. 
\end{proposition}

\textbf{Proof.} 
The relation 
\[g(t) = g_0(t) + \frac{1}{2(n-1)} \, t^2 \, u \, \overline{g}\] 
implies 
\[R_{g(t)} = R_{g_0(t)} - \frac{1}{2} \, t^2 \, (\Delta_{\overline{g}} u + nu) + O(t^3).\] 
Using Proposition \ref{key}, we obtain 
\[\frac{\partial}{\partial t} R_{g(t)} \Big |_{t=0} = 0\] 
and 
\[\frac{\partial^2}{\partial t^2} R_{g(t)} \Big |_{t=0} = Q - (\Delta_{\overline{g}} u + nu) = \mu\] 
at each point on $S_+^n$. This completes the proof. \\

Finally, we analyze the mean curvature of $\Sigma$ with respect to the metric $g(t)$.

\begin{proposition}
We have 
\[\frac{\partial}{\partial t} H_{g(t)} \Big |_{t=0} = -(\Delta_\Sigma \eta + (n-1)\eta) > 0\] 
at each point on $\Sigma$.
\end{proposition}

\textbf{Proof.} 
Note that 
\[g(0) = g_1(0) = \overline{g}\] 
and 
\[\frac{\partial}{\partial t} g(t) \Big |_{t=0} = \frac{\partial}{\partial t} g_1(t) \Big |_{t=0} = \mathscr{L}_X \overline{g}.\] 
This implies 
\[\frac{\partial}{\partial t} H_{g(t)} \Big |_{t=0} = \frac{\partial}{\partial t} H_{g_1(t)} \Big |_{t=0}\] 
at each point on $\Sigma$.

On the other hand, $H_{g_1(t)}$ can be identified with the mean curvature of the embedding $\varphi_t|_\Sigma: \Sigma \to S^n$ with respect to the standard metric $\overline{g}$. Moreover, we have $\frac{\partial}{\partial t} \varphi_t(x) \big |_{t=0} = \eta(x) \, \nu(x)$ for each point $x \in \Sigma$. Therefore, the standard formula for the linearization of the mean curvature gives 
\[H_{g_1(t)} = -t \, (\Delta_\Sigma \eta + (n-1)\eta) + O(t^2)\] 
at each point on $\Sigma$ (see e.g. \cite{Huisken-Polden}, Theorem 3.2). Putting these facts together, the assertion follows. \\

\begin{corollary}
If $t > 0$ is sufficiently small, then the scalar curvature of $g(t)$ is strictly greater than $n(n-1)$ at each point on $S_+^n$, and the mean curvature of $\Sigma$ with respect to $g(t)$ is strictly positive. Furthermore, we have $g(t) - \overline{g} = 0$ at each point on $\Sigma$.
\end{corollary}

\section{Proof of Theorem \ref{thm.B}}

Let $M$ be a compact Riemannian manifold of dimension $n$ with boundary $\partial M$. Suppose that $g$ and $\tilde{g}$ are two Riemannian metrics on $M$ with the property that $g - \tilde{g} = 0$ at each point on $\partial M$. Moreover, we assume that $H_g > H_{\tilde{g}}$ at each point on $\partial M$.

We will consider Riemannian metrics of the form $\hat{g} = g + h$, where $h$ is a suitably chosen perturbation. The following result provides an estimate for the scalar curvature of $\hat{g}$. In \cite{Brendle}, a similar result was established for perturbations of the Euclidean metric (cf. \cite{Brendle}, Proposition 26).

\begin{proposition}
\label{scalar.curvature}
Consider a Riemannian metric of the form $\hat{g} = g + h$, where $h$ satisfies the pointwise estimate $|h|_g \leq \frac{1}{2}$. Then the scalar curvature of $\hat{g}$ satisfies the estimate 
\begin{align*} 
&\bigg | R_{\hat{g}} - R_g - \sum_{i,j=1}^n (D_{e_i,e_j}^2 h)(e_i,e_j) + \Delta_g (\text{\rm tr}_g(h)) + \langle \text{\rm Ric}_g,h \rangle \bigg | \\ 
&\leq C \, |h|^2 + C \, |Dh|^2 + C \, |h| \, |D^2 h|. 
\end{align*} 
Here, $D$ denotes the Levi-Civita connection with respect to the metric $g$, and $C$ is a uniform constant which depends only on $(M,g)$. 
\end{proposition}

\textbf{Proof.} 
Let $\hat{D}$ be the Levi-Civita connection with respect to $\hat{g}$. Then 
\[\hat{D}_X Y = D_X Y + \Gamma(X,Y),\] 
where $\Gamma$ is defined by 
\begin{align*} 
2 \, \hat{g}(\Gamma(X,Y),Z) 
&= (D_X \hat{g})(Y,Z) + (D_Y \hat{g})(X,Z) - (D_Z \hat{g})(X,Y) \\ 
&= (D_X h)(Y,Z) + (D_Y h)(X,Z) - (D_Z h)(X,Y) 
\end{align*} 
(see e.g. \cite{Brendle-book}, Lemma A.2). In local coordinates, the tensor $\Gamma$ is given by 
\[\Gamma_{jk}^m = \frac{1}{2} \, \hat{g}^{lm} \, (D_j h_{kl} + D_k h_{jl} - D_l h_{jk}).\] 
The covariant derivatives of $\Gamma$ with respect to the metric $g$ are given by 
\begin{align*} 
D_i \Gamma_{jk}^m 
&= \frac{1}{2} \, \hat{g}^{lm} \, (D_{i,j}^2 h_{kl} + D_{i,k}^2 h_{jl} - D_{i,l}^2 h_{jk}) \\ 
&- \Gamma_{il}^m \, \Gamma_{jk}^l - \hat{g}^{lm} \, \hat{g}_{pq} \, \Gamma_{il}^q \, \Gamma_{jk}^p. 
\end{align*} 
The Riemann curvature tensor of $\hat{g}$ is related to the Riemann curvature tensor of $g$ by 
\begin{align*} 
(R_{\hat{g}})_{ijk}^m 
&= (R_g)_{ijk}^m + D_i \Gamma_{jk}^m - D_j \Gamma_{ik}^m + \Gamma_{jk}^l \, \Gamma_{il}^m - \Gamma_{ik}^l \, \Gamma_{jl}^m \\ 
&= (R_g)_{ijk}^m - \hat{g}^{lm} \, \hat{g}_{pq} \, \Gamma_{il}^q \, \Gamma_{jk}^p + \hat{g}^{lm} \, \hat{g}_{pq} \, \Gamma_{jl}^q \, \Gamma_{ik}^p \\ 
&+ \frac{1}{2} \, \hat{g}^{lm} \, (D_{i,j}^2 h_{kl} + D_{i,k}^2 h_{jl} - D_{i,l}^2 h_{jk}) \\ 
&- \frac{1}{2} \, \hat{g}^{lm} \, (D_{j,i}^2 h_{kl} + D_{j,k}^2 h_{il} - D_{j,l}^2 h_{ik}). 
\end{align*} 
Therefore, the scalar curvature of $\hat{g}$ is given by 
\begin{align*} 
R_{\hat{g}} 
&= \hat{g}^{ik} \, (\text{\rm Ric}_g)_{ik} 
+ \hat{g}^{ik} \, \hat{g}^{jl} \, \hat{g}_{pq} \, \Gamma_{il}^q \, \Gamma_{jk}^p - \hat{g}^{ik} \, \hat{g}^{jl} \, \hat{g}_{pq} \, \Gamma_{jl}^q \, \Gamma_{ik}^p \\ 
&- \frac{1}{2} \, \hat{g}^{ik} \, \hat{g}^{jl} \, (D_{i,j}^2 h_{kl} + D_{i,k}^2 h_{jl} - D_{i,l}^2 h_{jk}) \\ 
&+ \frac{1}{2} \, \hat{g}^{ik} \, \hat{g}^{jl} \, (D_{j,i}^2 h_{kl} + D_{j,k}^2 h_{il} - D_{j,l}^2 h_{ik}) \\ 
&= \hat{g}^{ik} \, (\text{\rm Ric}_g)_{ik} 
+ \hat{g}^{ik} \, \hat{g}^{jl} \, \hat{g}_{pq} \, \Gamma_{il}^q \, \Gamma_{jk}^p - \hat{g}^{ik} \, \hat{g}^{jl} \, \hat{g}_{pq} \, \Gamma_{jl}^q \, \Gamma_{ik}^p \\ 
&- \hat{g}^{ik} \, \hat{g}^{jl} \, (D_{i,k}^2 h_{jl} - D_{i,l}^2 h_{jk}). 
\end{align*} From this, the assertion follows easily. \\

We next describe our choice of perturbation. To that end, we fix a neighborhood $U$ of $\partial M$. Moreover, let us fix a smooth boundary defining function $\rho: M \to \mathbb{R}$ so that $\rho = 0$ and $|\nabla \rho| = 1$ at each point on $\partial M$. Since $g - \tilde{g}$ vanishes along $\partial M$, we can find a symmetric two-tensor $T$ such that $\tilde{g} = g + \rho \, T$ in a neighborhood of $\partial M$ and $T = 0$ outside $U$. The second fundamental form of $\partial M$ with respect to $\tilde{g}$ is given by  
\[A_{\tilde{g}}(X,Y) = A_g(X,Y) - \frac{1}{2} \, T(X,Y)\] 
for all vectors $X,Y \in T(\partial M)$. This implies 
\[H_{\tilde{g}} = H_g - \frac{1}{2} \, \text{\rm tr}(T|_{\partial M}).\] 
By assumption, we have $H_g > H_{\tilde{g}}$ at each point on $\partial M$. This implies $\text{\rm tr}(T|_{\partial M}) > 0$ at each point on $\partial M$. 

We next construct a suitable cut-off function:

\begin{lemma} 
\label{cutoff.function}
There exists a smooth cut-off function $\chi: [0,\infty) \to [0,1]$ with the following properties: 
\begin{itemize}
\item $\chi(s) = s - \frac{1}{2} \, s^2$ for each $s \in [0,\frac{1}{2}]$.
\item $\chi(s)$ is constant for $s \geq 1$.
\item $\chi''(s) < 0$ for all $s \in [0,1)$.
\end{itemize}
\end{lemma} 

\textbf{Proof.} 
We can find a smooth function $\psi: [0,\infty) \to \mathbb{R}$ such that $\psi(s) = 1$ for $s \in [0,\frac{1}{2}]$, $\psi(s) > 0$ for all $s \in [0,1)$, and $\psi(s) = 0$ for $s \geq 1$. Moreover, we may choose $\psi(s)$ such that $\int_0^\infty \psi(s) \, ds = 1$. We now define $\chi(s)$ as the unique solution of $\chi''(s) = -\psi(s)$ with initial conditions $\chi(0) = 0$ and $\chi'(0) = 1$. It is straightforward to verify that the function $\chi(s)$ has all the required properties. \\

Let $\beta: (-\infty,0] \to [0,1]$ be a smooth cutoff function such that $\beta(s) = \frac{1}{2}$ for $s \in [-1,0]$ and $\beta(s) = 0$ for $s \in (-\infty,-2]$. If $\lambda > 0$ is sufficiently large, we define a metric $\hat{g}_\lambda$ on $M$ by 
\[\hat{g}_\lambda = \begin{cases} 
g + \lambda^{-1} \, \chi(\lambda\rho) \, T & \text{\rm for $\rho \geq e^{-\lambda^2}$} \\ 
\tilde{g} - \lambda \rho^2 \, \beta(\lambda^{-2} \log \rho) \, T & \text{\rm for $\rho < e^{-\lambda^2}$}. \end{cases}\] 
If $\lambda > 0$ is sufficiently large, then $\hat{g}_\lambda$ is a smooth metric on $M$. Moreover, we have $\hat{g}_\lambda = \tilde{g}$ in the region $\{\rho \leq e^{-2\lambda^2}\}$ and $\hat{g}_\lambda = g$ outside $U$.

In the next step, we give a lower bound for the scalar curvature of $\hat{g}_\lambda$. We first consider the region $\{\rho \geq e^{-\lambda^2}\}$.

\begin{proposition}
\label{exterior.region}
Let $\varepsilon$ be an arbitrary positive real number. If $\lambda > 0$ is sufficiently large, then 
\[\inf_{\{\rho \geq e^{-\lambda^2}\}} (R_{\hat{g}_\lambda} - R_g) \geq -\varepsilon.\]
\end{proposition}

\textbf{Proof.} 
In the region $\{\rho \geq e^{-\lambda^2}\}$, we have $\hat{g}_\lambda = g + h_\lambda$, where 
\[h_\lambda = \lambda^{-1} \, \chi(\lambda \rho) \, T.\] 
The tensor $h_\lambda$ satisfies 
\begin{align*} 
\sum_{i,j=1}^n (D_{e_i,e_j}^2 h_\lambda)(e_i,e_j) 
&= \lambda \, \chi''(\lambda \rho) \, T(\nabla \rho,\nabla \rho) + \chi'(\lambda \rho) \, \langle D^2 \rho,T \rangle \\ 
&+ 2 \, \chi'(\lambda \rho) \, \sum_{j=1}^n (D_{e_j} T)(\nabla \rho,e_j) \\ 
&+ \lambda^{-1} \, \chi(\lambda \rho) \, \sum_{i,j=1}^n (D_{e_i,e_j}^2 T)(e_i,e_j) 
\end{align*} 
and 
\begin{align*} 
\Delta_g (\text{\rm tr}_g(h_\lambda)) 
&= \lambda \, \chi''(\lambda \rho) \, |\nabla \rho|^2 \, \text{\rm tr}_g(T) + \chi'(\lambda \rho) \, \Delta_g \rho \: \text{\rm tr}_g(T) \\ 
&+ 2 \, \chi'(\lambda \rho) \, \langle \nabla \rho,\nabla (\text{\rm tr}_g(T)) \rangle + \lambda^{-1} \, \chi(\lambda \rho) \, \Delta_g(\text{\rm tr}_g(T)). 
\end{align*}
Using Proposition \ref{scalar.curvature}, we obtain 
\begin{align*} 
&\Big | R_{\hat{g}_\lambda} - R_g + \lambda \, \chi''(\lambda \rho) \, \big ( |\nabla \rho|^2 \, \text{\rm tr}_g(T) - T(\nabla \rho,\nabla \rho) \big ) \Big | \\ 
&\leq N \, \lambda^{-1} \, \chi(\lambda \rho) + N \, \chi'(\lambda \rho) + N \, \chi(\lambda \rho) \, (-\chi''(\lambda \rho)) 
\end{align*} 
in the region $\{\rho \geq e^{-\lambda^2}\}$. Here, $N$ is a positive constant which is independent of $\lambda$. 

Recall that $\text{\rm tr}(T|_{\partial M}) > 0$ at each point on $\partial M$. By continuity, we can find a real number $a > 0$ such that 
\[|\nabla \rho|^2 \, \text{\rm tr}_g(T) - T(\nabla \rho,\nabla \rho) \geq a\] 
in a neighborhood of $\partial M$. Hence, if $\lambda > 0$ is sufficiently large, then we have  
\begin{align*} 
R_{\hat{g}_\lambda} - R_g 
&\geq -N \, \lambda^{-1} \, \chi(\lambda \rho) - N \, \chi'(\lambda \rho) + (a\lambda - N \, \chi(\lambda\rho)) \, (-\chi''(\lambda \rho)) \\ 
&\geq -N \, \lambda^{-1} - N \, \chi'(\lambda \rho) + (a\lambda - N) \, (-\chi''(\lambda \rho)) 
\end{align*}
in the region $\{\rho \geq e^{-\lambda^2}\}$. In the sequel, we always assume that $\lambda$ is chosen sufficiently large so that $a\lambda > N$.

Let us fix a real number $s_0 \in [0,1)$ such that $N \, \chi'(s_0) < \varepsilon$. Then 
\begin{align*} 
\inf_{\{e^{-\lambda^2} \leq \rho < s_0 \, \lambda^{-1}\}} (R_{\hat{g}_\lambda} - R_g) 
&\geq -N \, \lambda^{-1} - N + (a\lambda - N) \, \inf_{0 \leq s < s_0} (-\chi''(s)). 
\end{align*} 
By Lemma \ref{cutoff.function}, we have $\inf_{0 \leq s < s_0} (-\chi''(s)) > 0$. Thus, we conclude that 
\[\inf_{\{e^{-\lambda^2} \leq \rho < s_0 \, \lambda^{-1}\}} (R_{\hat{g}_\lambda} - R_g) \to \infty\] 
as $\lambda \to \infty$. Moreover, we have 
\[\inf_{\{\rho \geq s_0 \, \lambda^{-1}\}} (R_{\hat{g}_\lambda} - R_g) \geq -N \, \lambda^{-1} - N \, \sup_{s \geq s_0} \chi'(s) = -N \, \lambda^{-1} - N \, \chi'(s_0).\] 
Since $N \, \chi'(s_0) < \varepsilon$, it follows that 
\[\inf_{\{\rho \geq s_0 \, \lambda^{-1}\}} (R_{\hat{g}_\lambda} - R_g) \geq -\varepsilon\] 
if $\lambda > 0$ is sufficiently large. Putting these facts together, we conclude that 
\[\inf_{\{\rho \geq e^{-\lambda^2}\}} (R_{\hat{g}_\lambda} - R_g) \geq -\varepsilon\] 
if $\lambda > 0$ is sufficiently large. This completes the proof. \\

Finally, we estimate the scalar curvature of $\hat{g}_\lambda$ in the region $\{\rho < e^{-\lambda^2}\}$.

\begin{proposition}
\label{interior.region}
Let $\varepsilon$ be an arbitrary positive real number. If $\lambda > 0$ is sufficiently large, then  
\[\inf_{\{\rho < e^{-\lambda^2}\}} (R_{\hat{g}_\lambda} - R_{\tilde{g}}) \geq -\varepsilon.\]
\end{proposition}

\textbf{Proof.} 
In the region $\{\rho < e^{-\lambda^2}\}$, we have $\hat{g}_\lambda = \tilde{g} + \tilde{h}_\lambda$, where $\tilde{h}_\lambda$ is defined by 
\[\tilde{h}_\lambda = -\lambda\rho^2 \, \beta(\lambda^{-2} \log \rho) \, T.\] 
Let $\{e_1,\hdots,e_n\}$ denote a local orthonormal frame with respect to the metric $\tilde{g}$. Then we have 
\begin{align*} 
&\sum_{i,j=1}^n (\tilde{D}_{e_i,e_j}^2 \tilde{h}_\lambda)(e_i,e_j) \\ 
&= -\big [ 2\lambda \, \beta(\lambda^{-2} \log \rho) + 3\lambda^{-1} \, \beta'(\lambda^{-2} \log \rho) + \lambda^{-3} \, \beta''(\lambda^{-2} \log \rho) \big ] \, T(\tilde{\nabla} \rho,\tilde{\nabla} \rho) \\ 
&- \big [ 2\lambda\rho \, \beta(\lambda^{-2} \log \rho) + \lambda^{-1} \, \rho \, \beta'(\lambda^{-2} \log \rho) \big ] \, \langle \tilde{D}^2 \rho,T \rangle \\ 
&- \big [ 4\lambda\rho \, \beta(\lambda^{-2} \log \rho) + 2\lambda^{-1} \, \rho \, \beta'(\lambda^{-2} \log \rho) \big ] \, \sum_{j=1}^n (\tilde{D}_{e_j} T)(\tilde{\nabla} \rho,e_j) \\ 
&- \lambda\rho^2 \, \beta(\lambda^{-2} \log \rho) \, \sum_{i,j=1}^n (\tilde{D}_{e_i,e_j}^2 T)(e_i,e_j) 
\end{align*} 
and 
\begin{align*} 
&\Delta_{\tilde{g}} (\text{\rm tr}_{\tilde{g}}(\tilde{h}_\lambda)) \\ 
&= -\big [ 2\lambda \, \beta(\lambda^{-2} \log \rho) + 3\lambda^{-1} \, \beta'(\lambda^{-2} \log \rho) + \lambda^{-3} \, \beta''(\lambda^{-2} \log \rho) \big ] \, |\tilde{\nabla} \rho|^2 \, \text{\rm tr}_{\tilde{g}}(T) \\ 
&- \big [ 2\lambda\rho \, \beta(\lambda^{-2} \log \rho) + \lambda^{-1} \, \rho \, \beta'(\lambda^{-2} \log \rho) \big ] \, \Delta_{\tilde{g}} \rho \: \text{\rm tr}_{\tilde{g}}(T) \\ 
&- \big [ 4\lambda\rho \, \beta(\lambda^{-2} \log \rho) + 2\lambda^{-1} \, \rho \, \beta'(\lambda^{-2} \log \rho) \big ] \, \langle \tilde{\nabla} \rho,\tilde{\nabla} (\text{\rm tr}_{\tilde{g}}(T)) \rangle \\ 
&- \lambda\rho^2 \, \beta(\lambda^{-2} \log \rho) \, \Delta_{\tilde{g}}(\text{\rm tr}_{\tilde{g}}(T)). 
\end{align*} 
Using Proposition \ref{scalar.curvature}, we obtain
\[\Big | R_{\hat{g}_\lambda} - R_{\tilde{g}} - 2\lambda \, \beta(\lambda^{-2} \log \rho) \, \big ( |\tilde{\nabla} \rho|^2 \, \text{\rm tr}_{\tilde{g}}(T) - T(\tilde{\nabla} \rho,\tilde{\nabla} \rho) \big ) \Big | \leq L \, \lambda^{-1}\] 
in the region $\{\rho < e^{-\lambda^2}\}$. Here, $L$ is a positive constant which does not depend on $\lambda$.

Recall that $\text{\rm tr}(T|_{\partial M}) > 0$ at each point on $\partial M$. By continuity, we have 
\[|\tilde{\nabla} \rho|^2 \, \text{\rm tr}_{\tilde{g}}(T) - T(\tilde{\nabla} \rho,\tilde{\nabla} \rho) \geq 0\] 
in a neighborhood of $\partial M$. Hence, if $\lambda > 0$ is sufficiently large, then we have 
\[\inf_{\{\rho < e^{-\lambda^2}\}} (R_{\hat{g}_\lambda} - R_{\tilde{g}}) \geq -L \, \lambda^{-1}.\] From this the assertion follows. \\

Combining Proposition \ref{exterior.region} and Proposition \ref{interior.region}, we can draw the following conclusion: 

\begin{corollary}
Let $\varepsilon$ be a given positive real number. If we choose $\lambda > 0$ sufficiently large, then we have the pointwise inequality 
\[R_{\hat{g}_\lambda}(x) \geq \min \{R_g(x),R_{\tilde{g}}(x)\} - \varepsilon\] 
for each point $x \in M$.
\end{corollary}

\section{Proof of Theorem \ref{thm.C}}

In this final section, we describe the proof of Theorem \ref{thm.C}. As above, let $\overline{g}$ be the standard metric on $S^n$, and let $f$ denote the restriction of the coordinate function $x_{n+1}$.

\begin{lemma}
\label{subharmonic}
Assume that $\delta > 0$ is sufficiently small. Then the function $e^{-\frac{1}{f-\delta}}$ is subharmonic in the region $\{\delta < f < 3\delta\}$.
\end{lemma}

\textbf{Proof.} 
Assume that $0 < \delta < \frac{1}{8}$. Then $1 - 2(f-\delta) \geq \frac{1}{2}$ and $|\nabla f|^2 \geq \frac{1}{2}$ in the region $\{\delta < f < 3\delta\}$. This implies 
\begin{align*} 
\Delta_{\overline{g}} (e^{-\frac{1}{f-\delta}}) 
&= e^{-\frac{1}{f-\delta}} \, \bigg ( \frac{1 - 2(f-\delta)}{(f-\delta)^4} \, |\nabla f|^2 + \frac{1}{(f-\delta)^2} \, \Delta_{\overline{g}} f \bigg ) \\ 
&\geq e^{-\frac{1}{f-\delta}} \, \bigg ( \frac{1}{4 \, (f-\delta)^4} - \frac{nf}{(f-\delta)^2} \bigg ) 
\end{align*}
in the region $\{\delta < f < 3\delta\}$. Hence, if we choose $\delta > 0$ sufficiently small, then 
\[\Delta_{\overline{g}} (e^{-\frac{1}{f-\delta}}) \geq 0\] 
in the region $\{\delta < f < 3\delta\}$. \\

For each $\delta > 0$, we can find a smooth metric $\tilde{g}_\delta$ on $S_+^n$ such that 
\[\tilde{g}_\delta = \begin{cases} \overline{g} & \text{\rm for $f \leq \delta$} \\ (1 - e^{-\frac{1}{f-\delta}})^{\frac{4}{n-2}} \, \overline{g} & \text{\rm for $\delta < f < 3\delta$}. \end{cases}\] 
Using Lemma \ref{subharmonic}, we obtain a bound for the scalar curvature of the metric $\tilde{g}_\delta$.

\begin{proposition}
If $\delta > 0$ is sufficiently small, then the scalar curvature of $\tilde{g}_\delta$ is strictly greater than $n(n-1)$ in the region $\{\delta < f < 3\delta\}$.
\end{proposition}

\textbf{Proof.} 
Using the formula for the change of the scalar curvature under a conformal change of the metric, we obtain 
\begin{align*} 
R_{\tilde{g}_\delta} 
&= \frac{4(n-1)}{n-2} \, (1 - e^{-\frac{1}{f-\delta}})^{-\frac{n+2}{n-2}} \, \Delta_{\overline{g}}(e^{-\frac{1}{f-\delta}}) \\ 
&+ n(n-1) \, (1 - e^{-\frac{1}{f-\delta}})^{-\frac{4}{n-2}}
\end{align*}
in the region $\{\delta < f < 3\delta\}$. By Lemma \ref{subharmonic}, the function $e^{-\frac{1}{f-\delta}}$ is subharmonic in the region $\{\delta < f < 3\delta\}$. Thus, we conclude that 
\[R_{\tilde{g}_\delta} \geq n(n-1) \, (1 - e^{-\frac{1}{f-\delta}})^{-\frac{4}{n-2}} > n(n-1)\] 
in the region $\{\delta < f < 3\delta\}$. \\

For each $\tau \in (-1,1)$, we define a conformal diffeomorphism $\Psi_\tau: S^n \to S^n$ by 
\begin{align*} 
\Psi_\tau: \; &(x_1,\hdots,x_n,x_{n+1}) \\ 
&\mapsto \frac{1}{1+\tau^2+2\tau x_{n+1}} \, \Big ( (1-\tau^2)x_1, \hdots, (1-\tau^2)x_n,(1+\tau^2)x_{n+1}+2\tau \Big ). 
\end{align*}
If we choose $\tau = -\frac{2\delta}{1 + \sqrt{1-4\delta^2}}$, then $\Psi_\tau$ maps the domain $M_\delta = \{f \geq 2\delta\}$ to the hemisphere $S_+^n = \{f \geq 0\}$. 

By Theorem \ref{thm.A}, we can find a metric $g$ on $S_+^n$ with the following properties:
\begin{itemize}
\item The scalar curvature of $g$ is strictly greater than $n(n-1)$.
\item We have $g - \overline{g} = 0$ at each point on $\partial S_+^n$.
\item The mean curvature of $\partial S_+^n$ with respect to $g$ is strictly positive.
\end{itemize}
For $\delta > 0$ sufficiently small, we define a metric $g_\delta$ on $M_\delta$ by 
\[g_\delta = (1 - e^{-\frac{1}{\delta}})^{\frac{4}{n-2}} \, (1 - 4\delta^2) \, \Psi_\tau^*(g),\] 
where $\tau = -\frac{2\delta}{1 + \sqrt{1-4\delta^2}}$. Clearly, the scalar curvature of $g_\delta$ is strictly greater than $n(n-1)$. In the next step, we show that $g_\delta$ and $\tilde{g}_\delta$ agree along the boundary $\partial M_\delta$.

\begin{proposition}
\label{metrics.match}
We have $g_\delta - \tilde{g}_\delta = 0$ at each point on $\partial M_\delta$.
\end{proposition}

\textbf{Proof.} 
Using the relation $\delta = -\frac{\tau}{1+\tau^2}$, we obtain 
\[\Big ( \frac{1-\tau^2}{1+\tau^2+4\tau\delta} \Big )^2 = \Big ( \frac{1+\tau^2}{1-\tau^2} \Big )^2 = \frac{1}{1-4\delta^2}.\] 
It is straightforward to verify that 
\[\Psi_\tau^*(\overline{g}) = \Big ( \frac{1-\tau^2}{1+\tau^2+2\tau x_{n+1}} \Big )^2 \, \overline{g}\] 
at each point on $S^n$. This implies 
\[\Psi_\tau^*(\overline{g}) = \Big ( \frac{1-\tau^2}{1+\tau^2+4\tau\delta} \Big )^2 \, \overline{g} = \frac{1}{1-4\delta^2} \, \overline{g}\] 
at each point on $\partial M_\delta$. Since $g$ agrees with $\overline{g}$ along the boundary $\partial S_+^n$, we conclude that 
\[g_\delta = (1 - e^{-\frac{1}{\delta}})^{\frac{4}{n-2}} \, (1 - 4\delta^2) \, \Psi_\tau^*(\overline{g}) = (1 - e^{-\frac{1}{\delta}})^{\frac{4}{n-2}} \, \overline{g} = \tilde{g}_\delta\] 
at each point on $\partial M_\delta$. This completes the proof of Proposition \ref{metrics.match}. \\

To conclude the proof of Theorem \ref{thm.C}, we choose $\delta > 0$ sufficiently small so that 
\[\sup_{\partial M_\delta} H_{\tilde{g}_\delta} < \inf_{\partial S_+^n} H_g.\] 
This implies 
\[\sup_{\partial M_\delta} H_{\tilde{g}_\delta} < \inf_{\partial M_\delta} H_{g_\delta}.\] 
By Theorem \ref{thm.B}, there exists a smooth Riemannian metric $\hat{g}$ on $M_\delta$ which has scalar curvature strictly greater than $n(n-1)$ and which agrees with $\tilde{g}_\delta$ in a neighborhood of $\partial M_\delta$. Hence, the metric $\hat{g}$ extends to a smooth metric on the hemisphere with the property that $R_{\hat{g}} \geq n(n-1)$ at each point on $S_+^n$ and $\hat{g} = \tilde{g}_\delta$ in the region $\{f \leq 2\delta\}$. In particular, we have $\hat{g} = \overline{g}$ in the region $\{f \leq \delta\}$. This completes the proof of Theorem \ref{thm.C}.

\end{document}